\documentclass[10pt,a4paper]{amsart}
\usepackage{amssymb}
\usepackage{amscd}
\usepackage{amsmath,amsfonts,stmaryrd}
\input xy 
\xyoption{all}

\usepackage{lscape}
\usepackage[colorlinks=true,linkcolor=black,pdftex]{hyperref}
\usepackage[hyperpageref]{backref} 
\usepackage{faktor}
\usepackage[nice]{nicefrac}
\usepackage{calc}
\usepackage{niceframe}
\theoremstyle{plain}

\newtheorem{theorem}{Theorem}[section]

\newtheorem{lemma}{Lemma}
\newtheorem{proposition}{Proposition}
\theoremstyle{definition}
\newtheorem{definition}{Definition}
\theoremstyle{definition}

\newtheorem*{remark}{Remark}

\theoremstyle{remark}

\DeclareMathOperator{\id}{id}
\DeclareMathOperator{\Inc}{Inc}
\DeclareMathOperator{\Av}{Av}

\DeclareMathOperator{\EXT}{\rm Ext}

\DeclareMathOperator{\Hom}{\rm Hom}
\DeclareMathOperator{\End}{\rm End}

\DeclareMathOperator{\rk}{\rm rk}
\DeclareMathOperator{\inj}{\hookrightarrow}

\newcommand{\oTo}{\xymatrix{ \ar@{^{(}->}[r]|{\mathbf{O}}& }}
\newcommand{\cTo}{\xymatrix{ \ar@{^{(}->}[r]|{\mathbf{|}}& }}
\newcommand{\coTo}{\xymatrix{ \ar@{^{(}->}[r]|{\mathbf{O}}|{\mathbf{|}}& }}
\DeclareMathOperator{\surj}{\twoheadrightarrow }

\DeclareMathOperator{\Bild}{Im}

\DeclareMathOperator{\Lie}{Lie}

\newcommand{\Gl}{\mathbf{Gl}}

\DeclareMathOperator{\Grmod}{grmod}

\newcommand{\si}{\sigma}

\newcommand{\al}{\alpha}

\newcommand{\B}{\mathbb{B}}
\newcommand{\C}{\mathbb{C}}

\newcommand{\G}{\mathbb{G}}

\newcommand{\N}{\mathbb{N}}

\renewcommand{\P}{\mathbb{P}}

\newcommand{\U}{\mathbb{U}}

\newcommand{\W}{\mathbb{W}}

\newcommand{\Z}{\mathbb{Z}}

\newcommand{\mcC}{\mathcal{C}}

\newcommand{\mcF}{\mathcal{F}}
\newcommand{\mcG}{\mathcal{G}}

\newcommand{\mcL}{\mathcal{L}}

\newcommand{\mcN}{\mathcal{N}}

\newcommand{\mcU}{\mathcal{U}}

\newcommand{\mcZ}{\mathcal{Z}}

\newcommand{\NH}{\mathrm{NH}}

\newcommand{\ddd}{\underline{\mathbf{d}}}

\newcommand{\Gd}{\mathbf{Gl}_{\underline{d }}}

\begin{document}

\title{From complete to partial flags in geometric Extension algebras}
\author{Julia Sauter} 
\address{Faculty of Mathematics, Bielefeld University}
\email{jsauter@math.uni-bielefeld.de}



\begin{abstract}
A geometric extension algebra is an extension algebra of a semi-simple perverse sheaf (allowing shifts), e.g. a push-forward of the constant sheaf under a projective map. Particular nice situations arise for collapsings of homogeneous vector bundle over homogeneous spaces. In this paper, we study the relationship between partial flag and complete flag cases. Our main result is that the locally finite modules over the geometric extension algebras are related by a recollement. As examples, we investigate parabolic affine nil Hecke algebras, geometric extension algebras associated to parabolic Springer maps and an example of Reineke of a parabolic quiver-graded Hecke algebra.   
\end{abstract}
\maketitle

\section{Introduction}
A geometric extension algebra is an algebra $\EXT^*(\mcL, \mcL ) := \bigoplus_{n\in \Z} \Hom (\mcL, \mcL [n])$ where $\mcL$ is a direct sum of shifts of simple equivariant 
perverse sheaves on a complex algebraic variety equipped with a suitable group action.  For example $\mcL$ might be a push-forward of the constant sheaf under an equivariant projective map of varieties. 
Examples of such algebras include affine nil Hecke algebras, skew group algebras of polynomial rings with a  Weyl group action, and graded parts of quiver Hecke algebras, see for example \cite{Sa2}. All of these algebras arise from collapsing of homogeneous vector bundles over complete flag varieties, and in this case it is often possible to find explicit generators and relations using methods of Varagnolo and Vasserot,  see for example \cite{VV}, \cite{VV2} and \cite{Sa}.
However, there are also cases of interest associated to partial flag varieties, to which these methods do not apply. This includes quiver Schur algebras in \cite{SW}, where calculations can be made via diagram calculus.  

In this paper we study geometric extension algebras $\mcZ^P$ arising from partial flag varieties. Associated to such an algebra we construct another 
geometric extension algebra $\mcZ^B$ arising from complete flags and 
we show that there is a close relation between these two algebras. 
More precisely, we show there is an idempotent element $e_P\in \mcZ^B$ such that $e_P\mcZ^Be_P \cong \mcZ^P$. This implies that the locally finite modules over $\mcZ^B$ and $\mcZ^P$ are related by a recollement (obtained from the idempotent $e_P$). Our result applies to arbitrary reductive groups and more general collapsings than studied in loc.\ cit. As examples, we study parabolic versions of affine nil Hecke algebras, examples from Springer theory and an example of Reineke.

\section{A recollement relating parabolic and Borel cases.}

\begin{definition}\label{d:setup}
Let $(G,P_i, V, F_i)_{i\in I}$ be a tuple with $G$ a reductive group with parabolic subgroups $P_i, i\in I$ (where $I$ is some finite set) such 
that $\bigcap P_i$ contain a maximal torus $T$, $\dim P_i =\dim P_j$ for all $i,j\in I$ and $V$ a $G$-representation with $P_i$-subrepresentations $F_i \subset V$. 
We set 
$E_i^P:= G\times^{ P_i}F_i$ ,
$E^P:= \bigsqcup E_i^P , \pi^P\colon E^P\to V,\overline{(g,f)}\mapsto gf$. 
Choose $T\subset B_i\subset P_i,i\in I$ Borel subgroups of $G$ 
(where $T\subset \bigcap_{i\in I}P_i$) and consider $F_i$ as $B_i$-representation, then $(G,B_i,V,F_i)_{i\in I}$ can be used 
to define $E_i^B, E^B, \pi^B, Z^B$ analogously.
\end{definition}

For any complex algebraic variety $X$ with an action of an algebraic group $A$ we will denote by $D^b_A(X)$ the $A$-equivariant derived category introduced by 
Bernstein and Lunts \cite{BL}. These categories carry a six functor calculus, we will denote the right or left derived functors with the same symbol as the functor. They have a dualizing object $\textbf{D}_X$ used to define a duality called Verdier duality. 

We define $\EXT^*_{D^b_A(X)}(U,V):= \bigoplus_{n\in \Z} \Hom_{D_A^b(X)}(U,V[n])$ for $U,V\in D_G^b(X)$, and  $\mcC^B:= \bigoplus_{i\in I} \C_{E^B_i}[\dim E_i^B]$ and $\mcC^P$ analogously and we set 
\[
\begin{aligned}
 \mcZ^B &=\EXT^*_{D^A_b(V)}(\pi^B_*\mcC^B, \pi^B_*\mcC^B) \\
\mcZ^P &=\EXT^*_{D^A_b(V)}(\pi^P_*\mcC^P, \pi^P_*\mcC^P) 
\end{aligned}
\]
We remark that the shifts in the grading of the constant sheaves ensures $\textbf{D}_{E^B} (\mcC^B )=\mcC^B$ and 
$\textbf{D}_{E^P} (\mcC^P ) =\mcC^P$. This implies that both graded algebras $\mcZ^B, \mcZ^P$ have an anti-involution given by applying Verdier duality.  
We denote by $*-\Grmod$ the category of finitely generated graded modules over a $\Z$-graded algebra $*$. Homomorphisms are given by homogeneous maps of degree zero. 
Since $\mcZ^B$ and $\mcZ^P$ are finitely generated modules over the commutative noetherian ring $H_A^ *(pt)$ (see \cite{Sa2}), they are noetherian and the catgeories of graded modules over them are abelian. 

\begin{theorem}\label{t:main}
There is a recollement of abelian categories given by an idempotent element $e_P$ 
\[
\xymatrix{
\mcZ^B_A/\mcZ^B_A e_P\mcZ^B_A\!-\!\Grmod \ar[rr] \;  && \; \mcZ^B_A\!-\!\Grmod  \; \ar@<2ex>[ll] \ar@<-2ex>[ll] \ar[rr]|{e_P }  && \; \mcZ^P_A\!-\!\Grmod \ar@<2ex>[ll] \ar@<-2ex>[ll]
}
\]
\end{theorem}

In fact, given any graded algebra $B$ and an idempotent $e\in B$ in degree $0$, 
such that $B$ is a module-finite algebra over a commutative noetherian ring, we get that the standard recollement for the idempotent $e$ restricts to one of the finitely generated graded modules (using if $X$ and $Y$ are finitely generated $Z$-modules, then so is $\Hom_Z (X,Y)$). 

The key tool to find the idempotent $e_P$ in the theorem is a (topological) operation of the Weyl group.

\subsection{They Weyl group operation}
Let $T\subset B\subset P\subset G$ with $T$ be a maximal torus and $B$ be a Borel in a reductive group $G$ over the $\C$. We choose a maximal compact subgroup $K\subset G$, then $T^\prime:=T\cap K$ is a compact torus in $K$ and the inclusion $K\to G$ induces a 
homoemorphism $K/T^\prime \to G/B$, this is long well-known see e.g. \cite[section 2]{Ar3}. The Weyl group $W$ associated to $G$ and $T$ coincides with the Weyl group associated to $K$ and $T^\prime$. The group $W$ operates on $K/T^\prime$ via $nT^\prime \cdot kT^\prime = kn^{-1}T^\prime , n\in N_{K}(T^\prime ), k\in K$ without fixpoints. 
Similarly, let $W_P\subset W$ be the Weyl group of a Levi subgroup $L\subset P$ and $T$, 
it operates on $P/B \cong L/ (B\cap L) \cong K^\prime /T^\prime$ for a maximal compact subgroup $K^\prime \subset L$. Now given, a $P$-representation $F$, 
then $E^B:=G\times^B F\cong K\times^{T^\prime} F$ carries a $W$-operation and for the 
map $\al \colon E^B \to E^P:=G\times^PF, \overline{(g,f)}^B\mapsto \overline{(g,f)}^P$ one has $\al\circ w =\al$ for all $w\in W_P$. This induces an operation of $W_P$ on $\al_*(\mcF)$ for any sheaf $\mcF$ and therefore we have an induced functor 
\[ \al_* ()^{W_P} \colon D^b_A(E^B) \to D^b_A(E^P). \]

\begin{lemma}\label{keylemma} In the notation from before, for the constant sheaf $\underline{\C}$ on $E^P$ the adjunction map $\underline{\C}\to \al_*\al^*\underline{\C}$ is a monomorphism in $D_A^b(E^P)$. 
Furthermore it factorizes over an isomorphism $\underline{\C} \to (\al_* \underline{\C})^{W_P}$ where $W_P$ is the Weyl group of a Levi subgroup in $P$.  
\end{lemma}

\paragraph{proof:}
For any variety we write $X_A:= X\times^AEA$ where $EA$ is a contractible space with a free $A$-operation.  
We denote $\al:=\al_A\colon (E^B)_A\to (E^P)_A$ the associated map. It is smooth and proper submersion with fibres all isomorphic to $P/B$.  
We prove the lemma first for the constant sheaf $\mcF=\underline{\C}$. 
The decomposition theorem (in the more specific version for a proper submersion, see \cite{CM},p.14 3rd Example) implies
\[ \al_*\underline{\C} =\bigoplus_{i \in \Z} R^i\al_*\underline{\C}[-i]. \]
Since $R^i\al_*\underline{\C}$ is the sheaf associated to the presheaf 
\[
U\mapsto H^i(\al^{-1}(U))
\]
this implies $(R^i\al_*\underline{\C})_x=H^i(\al^{-1}(x)) \cong H^i(P/B)$ for all $x\in (E^P)_A$. 
Therefore, $R^i\al_*\underline{\C}$ is a local system on $(E^P)_A$ and since $\pi_1 ((E^P)_A, x_0)$ (for any $x_0\in (E^P)_A$) is trivial, 
it is the constant local system. Then, one has $  \al_*\underline{\C}= \bigoplus_{w\in W_P} \underline{\C} [-2 \ell (w)]$ because $H^*(P/B)= H^*(L/(L\cap B)=\C[ \mathfrak{t}]/I_{W_P}$ where the last 
isomorphism is graded algebras and as $W_P$-representations by the Borel isomorphism. 
But since $(\C[ \mathfrak{t}]/I_{W_P})^{W_P}=\C$ in degree $0$,  
one has $(\al_*\underline{\C})^{W_P}\cong \underline{\C}$. Furthermore, it is easy to see that the unit of the adjunction is a monomorphism 
(since $\al$ is locally trivial). By taking the trivial $W_P$-operation on $\underline{\C}$, we can make the unit of the adjunction a $W_P$-linear map (because the map is locally trivial and $W_P$ operates only on the fibre), then taking $W_P$-invariants proves the lemma. 
\hfill $\Box $

Now, let us come back to the setup from definition \ref{d:setup}. 
Let $i\in I$. Consider the following commutative triangle  
\[ 
\xymatrix{
E_i^B=G\times^{B_i}F_i \ar[rr]^{\al_i}\ar[rd]_{\pi^B_i} &&  E_i^P=G\times^{P_i}F_i \ar[ld]^{\pi^P_i}\\
& V& }
\]
where $\pi_i^B=\pi^B|_{E_i^B}$, $\pi_i^P=\pi^P|_{E_i^P}$. We denote by $W_i$ the Weyl group of a Levi subgroup in $P_i$ and the torus $T$. 
We define two morphisms in $D^b_A(V)$ using the isomorphisms from lemma \ref{keylemma}. 
\[
\begin{aligned}
\Inc_i \colon (\pi^P_i)_*\underline{\C}_{E^P_i} \cong (\pi_i^P)_*[(\al_i)_*\al_i^*\underline{\C}_{E_i^P}^{W_i}]
&\inj (\pi_i^P)_*(\al_i)_*\al_i^*\underline{\C}_{E_i^P} = 
(\pi^B_i)_*\underline{\C}_{E^B_i}\\
\Av_i\colon (\pi^B_i)_*\underline{\C}_{E^B_i}= 
(\pi_i^P)_*(\al_i)_*\al_i^*\underline{\C}_{E_i^P} 
& \surj (\pi_i^P)_*[(\al_i)_*\al_i^*\underline{\C}_{E_i^P}^{W_i}]= 
(\pi^P_i)_*\underline{\C}_{E^P_i} 
\end{aligned}
\]
The second map is locally the averaging map (or Reynolds operator) for the finite group $W_i$.  
One has $\Av_i\circ \Inc_i = \id$ and $\Inc_i\circ \Av_i =: e_i$ is an idempotent endomorphism.
Now let $s:= \dim P_i -\dim B_i$ for one (and all) $i\in I$ and observe $ \pi_*^P \mcC^P [s] = \bigoplus_{i\in I} (\pi_i^P)_* \underline{\C}_{E_i^P}[\dim E_i^B]$, we define 
\[
\begin{aligned}
\Inc:= \bigoplus_{i\in I} \Inc_i [\dim E_i^B] \colon \pi_*^P \mcC^P [s]& \to \pi_*^B\mcC^B \\
\Av:= \bigoplus_{i\in I} \Av_i[\dim E_i^B] \colon \pi_*^B \mcC^B &\to \pi_*^P\mcC^P[s]  
\end{aligned}
\]
and we observe that $\Av\circ \Inc= \id$ and $\Inc \circ \Av =e_P$ with   
 \[e_P:= (e_i)_{i\in I}\in \bigoplus_{i\in I}\End_{D^b_A(V)}( (\pi^B_i)_*\underline{\C}) = \bigoplus_{i\in I}\End_{D^b_A(V)}( (\pi^B_i)_*\underline{\C}[\dim E_i^B]) \subset (\mcZ^B)_0 .\]
Then, we define a map  
\[
 \theta_B^P \colon \mcZ^P= \EXT^*(\pi_*^P\mcC^P[s],\pi_*^P\mcC^P[s])  \to \mcZ^B 
 \]
 mapping an element $f$ in degree $n$ to $\Inc [n] \circ f\circ \Av$. It is easy to see that this map preserves degrees, $\theta_B^P(1) =e_P$ and products are mapped to products since $\Av \circ \Inc =\id$.

\begin{lemma}
The map $\theta_B^P$ is injective and induces an isomorphism of graded algebras 
\[ 
\mcZ^P \cong \Bild (\theta_B^P) =e_P \mcZ^B e_P. 
\] 
\end{lemma}

\paragraph{proof:} Assume $\Inc [n] \circ f \circ \Av =0$, then since $\Inc [n]$ is a monomorphism 
$f\circ \Av=0$. But then $f = f\circ (\Av \circ \Inc) =0$. The rest is clear.  
\hfill $\Box$

Now, the lower recollement is the standard recollement induced by the idempotent element. The map $e_P$ is left multiplication with the idempotent $e_P$. This completes the proof of theorem \ref{t:main}. 

Let $W_P:=\prod_{i\in I} W_i$, then $W_P$ operates on $\pi_*^B \mcC^B$ with $(\pi_*^B \mcC^B)^{W_P}=\pi_*^P\mcC^P [s]$, see before. 
This induces a $W_P\times W_P$-operation on $\mcZ^B$, for  
$f$ of degree $n$, $v,w\in W_P$ we set 
\[ (v,w) \cdot f := v[n] \circ f \circ w^{-1} .\]
This operation is $H_A^*(pt)$-linear but does not preserve 
products. 

\begin{proposition}\label{p:woperation}
$\theta_B^P(\mcZ^P) = (\mcZ^B)^{W_P\times W_P}$. 
\end{proposition}

\begin{proof}
Let $f\in \mcZ^P$ of degree $n$, then $\theta_B^P(f)= \Inc[n]\circ f\circ\Av$ is $W_P\times W_P$-invariant,  therefore $\Theta_B^P(\mcZ^P)\subset (\mcZ^B)^{W_P\times W_P}$.  
On the other hand, given $h\in \mcZ^B$ of degree $n$ with $ v[n]\circ h\circ w^{-1} = h, v,w\in W_P$, then one has $h\circ (\Inc \circ \Av) = (\Inc \circ \Av)[n] \circ h$ which implies 
$h ((\pi_*^B\mcC^B)^{W_P})\subset (\pi_*^B \mcC^B)^{W_P}[n] $, therefore restriction induces an element $\overline{h} \in \mcZ^P$ and by definition $h = \Inc\circ \overline{h}\circ \Av$.  
\end{proof}


\section{The hypercohomology functor}
Let $X$ be a complex variety with the action of an algebraic group $A$. We set\\
 $H^*_A(X):=H^*_A(X,\C) :=  \EXT^*_{D^b_A(X)}(\underline{\C}_X, \underline{\C}_X)$ for $A$-\textbf{equivariant cohomology} and \\
$H_*^A(X):= H_*^A(X,\C):= \EXT^{-*}_{D^b_A(X)}(\underline{\C}_X, \textbf{D}_X)$ for $A$-\textbf{equivariant Borel-Moore homology} \\
with complex coefficients respectively.

We associate to the data from before the Steinberg variety 
\[ Z^P:= \bigsqcup_{i,j\in I} \underbrace{E_i^P \times_V E_j^P}_{Z_{i,j}}  \text{ where } E_i^P := (G\times^{P_i}F_i), i\in I \]
and call $H_*^A(Z^P), $ $A\in \{pt, T, G\}$ the Steinberg algebra associated to the data (where the product is given by a convolution construction defined by \cite{CG}, section 2.7). It is a graded algebra with respect to 
\[
H_{[p]}^A(Z^P):= \bigoplus_{i,j\in I} H_{r_i+r_j-p}^A(Z_{i,j}^P) , \text{ where } r_i:= \dim_{\C} E_i^P , \]
we write $H_{[*]}^A(Z)$ to indicate this  grading. We recall the following result. 
\begin{theorem} \label{Extalg} (\cite[chapter 8]{CG})
 Let $A\in \{pt, T,G\}$ we write $\pi_i\colon E_i^P\to V, \overline{(g,f)}\mapsto gf$ and 
 There is an isomorphism of graded  $\C$-algebras 
\[ H_{[*]}^A(Z^P) \to \mcZ^P \]
\end{theorem}

For every $\mcF$ in $D^b_A(E^P)$, the hypercohomology $\mathbb{H}_A^*(E^P, \mcF)$ is naturally a bimodule over the equivariant cohomology ring $H_A^*(E^P)$. This is (a generalization of) S\"orgel's bifunctor which leads to the definition of 
S\"orgel bimodules. We think that the (lesser known) operation of the Steinberg algebra should also be of interest. 

\begin{theorem} For every $\mcF \in D^b_A(E^P)$ its hypercohomology $\mathbb{H}_A^*(E^P, \mcF)$ is a graded left and right module over $\mcZ^P$ with finite-dimensional graded parts. This module structure is natural in $\mcF$.  
\end{theorem}

The (ungraded) version of the convolution operation on hypercohomology groups of objects in the equivariant derived category can adapted from the main result of \cite{Jos}, the grading follows from the definition of the convolution operation, we explain it shortly. 
We shorten the notation here $E^P=:E, Z^P=:Z$. Recall $Z=E\times_V E$, set $E=E\times pt$ . 
Consider the $A$-equivariant maps 
\[
\xymatrix{ & E_i \times E_j\times pt \ar[ld]_{p_{12}} \ar[d]^{p_{23}} \ar[dr]^{p_{13}} & \\ 
E_i\times E_j & E_j\times pt & E_i\times pt 
}
\]
then since $\iota\colon Z_{ij}\subset E_i\times E_j$ is a closed embedding one has $\mathbb{H}_A^{*}(Z_{i,j}, D_{Z_{ij}})=
\mathbb{H}_A^{*}(E_i\times E_j, \iota _*(D))$ where $D$ is the dualizing sheaf and we define 
\[
\mathbb{H}_A^{p-r_i-r_j}(Z_{i,j}, D) \otimes \mathbb{H}^{k+r_t}_A(E_t,\mcF|_{E_t}) \xrightarrow{\delta_{j,t}(p_{13})_*(p_{12}^*[-]\otimes p_{23}^*[-]) } \mathbb{H}^{k+r_i}_A(E_i, \mcF|_{E_i})
\]
using the operators explained in \cite{Jos} section 2 and section 4.
The proof of the main result in \cite{Jos} adapts straight forward to this situation and proving that this map provides the left operation (for the right operation permute the factors). If we set 
$\mathbb{H}^{[t]}_A(E,\mcF) := \bigoplus_{i\in I}  \mathbb{H}^{t+r_i}_A(E_i,\mcF|_{E_i})$ this gives it the structure of a graded module. 

For $P_i=B_i$ for all $i\in I$, we write $Z^B, E^B,..$ instead of $Z^P, E^P, ..$. 
We observe that $W_P\times W_P$ operates on $Z^B$ by homeomorphisms and this corresponds under the automorphism in theorem \ref{Extalg} to the $W_P\times W_P$-operation from the previous section. Furthermore, if we restrict the operation to $W_P\times \{1\}$ and consider the $W_P$-operation on $E^B$, then the multiplication above commutes with this $W_P$-operation
on $\mathbb{H}^{k+r_i}_A(E_i^B, \mcF|_{E_i})= \mathbb{H}^{k+r_i}_A(E_i^P, (\al_i)_*\mcF|_{E_i})$
 i.e. for $x \in  \mathbb{H}_A^{p-r_i-r_j}(Z_{i,j}, D), f\in \mathbb{H}^{k+r_t}_A(E_t,\mcF|_{E_t}), w \in W_i$ we have 
\[
w (xf) = ((w,1) x) f \in \mathbb{H}^{k+r_i}_A(E_i, \mcF|_{E_i}),
\]
and in particular, $(\mcZ^B)^{W_P\times \{1\}}) \mathbb{H}^*_A(E^B, \mcF) = \mathbb{H}^*_A(E^P, \al_*(\mcF))^{W_P}$. 

So taking hypercohomology provides us with functors 
\[ 
\begin{aligned}
\mathbb{H}^*\colon D^b_A(E^B) & \to \mcZ^B\!-\!\Grmod, \\
\mathbb{H}^*\colon D^b_A(E^P) & \to \mcZ^P\!-\!\Grmod, 
\end{aligned}
\]
The complex $\underline{\C}$ maps to the equivariant 
comohomology $H^{[*]}_A(E^B)$ and $H^{[*]}_A(E^B)$. 

If one has $F_i^T=\{0\}$ for every $i\in I$ and $A\in\{T, G\}$ then $H^{[*]}_A(E^P)$ and $H^{[*]}_A(E^B)$ are faithful and for $A=pt$ this is not the case, see \cite{Sa2}.

So, combining our previous result with the hypercohomology functor, we can prove 
\begin{proposition} There is a commutative diagram 
\[
\xymatrix{ D^b_A( E^B)\ar[r]^{\al_*()^{W_P}}\ar[d]^{\mathbb{H}^*} & D^b_A(E^P)\ar[d]^{\mathbb{H}^*} \\
\mcZ^B\!-\!\Grmod \ar[r]^{e_P}& \mcZ^B\!-\!\Grmod
}
\]
\end{proposition}

\begin{proof}
We have by our previous observation  
$\mathbb{H}^*(\al_*(\mcF)^{W_P}) = \mathbb{H}^* (\al_*(\mcF))^{W_P}= 
(\mcZ^B)^{W_P\times \{1\}}) \mathbb{H}^*_A( \mcF)=e_P\cdot \mathbb{H}^*(\mcF)$
where the last equality is proven in the proof of Proposition \ref{p:woperation}. 
\end{proof}

%

\section{Examples and applications}
 
\subsection{The parabolic affine nil Hecke algebra}
Let $G$ be a reductive group over $\C$ and $B\subset G$ be a Borel subgroup. 
The \textbf{affine nil Hecke algebra} is defined by the graded vector space  
\[\NH:= \bigoplus_{p \in \Z} H_{[p]}^G(G/B\times G/B) \] 
where $H_{[p]}^G (G/B \times G/B):= H_{2\dim (G/B) -p}^G(G/B \times G/B)$, then $\NH$ is a graded algebra with respect to the convolution product defined by Chriss and Ginzburg (see \cite{CG}).   
We define  
\[\NH^P:=\bigoplus_{p \in \Z} H_{[p]}^G(G/P\times G/P)\]
with $H_{[p]}^G(G/P \times G/P ) := H_{2 \dim (G/P) -p}^G(G/P \times G/P)$, again this is a graded algebra with respect to the convolution product defined by Chriss and Ginzburg. We call it the \textbf{parabolic affine nil Hecke algebra}, it is a graded $H_G^*(pt)$-algebra. The following lemma has been observed in loc. cit. in the not equivariant case. 

\begin{lemma} 
One has $\NH^P \cong \End_{H_G^*(pt)}(H_G^*(G/P))$ as $\Z$-graded $H_G^*(pt)$-algebras.  
\end{lemma}

Recall that $H_G^*(pt) = (H_T^*(pt))^W= \C[\mathfrak{t}]^W$ is a commutative and graded $\C$-algebra 
where $T\subset P$ is a maximal torus, 
$\mathfrak{t}$ its Lie algebra and $W$ the Weyl group for $(G,T)$. 
Also we know that $H_G^*(G/P)\cong \C[\mathfrak{t}]^{W_P}$ 
where $W_P$ is the Weyl group of $(L,T)$ for the Levi subgroup $L\subset P$. We write $W^P\subset W$ for the minimal coset representatives of the cosets $W / W_P$. 

We give a proof of the previous lemma on the grounds that we could not find one in the literature. 
\begin{proof}
Let $EG$ be a contractible free $G$-space (or an appropriate approximation of it in the sense of \cite{BL}). 
Let $X:= G/P, \pi\colon X_G:= X\times^G EG \to BG $ the map obtained from $X\to pt$ by applying $-\times^G EG$. By \cite{CG}, chapter 8, 
we know $H_{[*]}^G(G/P\times G/P)\cong \EXT_{D^b_G(pt)}^* (\pi_* \underline{\C},\pi_* \underline{\C})$ as graded $H_G^*(pt)$-algebras. 
Since $\pi$ is a proper submersion, we have (by \cite{CM}, p.14 3rd Example) 
\[  \pi_*\underline{\C} = \bigoplus_{i\in \Z} R^i \pi_*\underline{\C}[-i] \]
in $D_G^b(pt)$. Since all fibres of $\pi$ are isomorphic to $X$ and $BG$ is simply connected, we get that 
\[\pi_*\underline{\C} =\bigoplus_{w\in W^P} \underline{\C}[-2\ell (w)],\] 
where $\ell (w)$ is the length of $w$. Let $r=\dim_{\C} H^*(X) =\# W^P$. 
We know, that $\C[\mathfrak{t}]^{W_P}$ is a free module over $\C[\mathfrak{t}]^W$ of rank $r$ generated by elements $b_w, w\in W^P, deg b_w = 2\ell (w)$. 
Now, the claim follows from the (wellknown) algebra isomorphism  
$\EXT_{D^b_G(pt)}^* (\underline{\C}, \underline{\C}) \cong H_G^*(pt)$. 
\end{proof}

Our maps from the earlier section give a natural homomorphism of graded $H_G^*(pt)$-modules
\[
\begin{aligned}
\Theta\colon \NH^P =\End_{H_G^*(pt)}(H_G^*(G/P))&\to \NH=\End_{H_G^*(pt)}(H_G^*(G/B))\\
f &\mapsto \Inc \circ f \circ \Av
\end{aligned}
\]  
where $\Inc \colon \C[\mathfrak{t}]^{W_P}\subset \C[\mathfrak{t}]$ is the natural inclusion and $\Av\colon \C[\mathfrak{t}]\to \C[\mathfrak{t}]^{W_P}, f\mapsto 
\frac{1}{\# W_P} \sum_{w\in W_P} w(f)$ 
is the averaging map. We set $e_P=\Inc\circ \Av (1)$.  
Furthermore, $W_P\times W_P$ operates on $\NH$ via graded $H_G^*(pt)$-module homomorphisms defined by $(v,w)\cdot h (f) := v(h(w^{-1}(f))), \quad v,w\in W_P, h\in \NH, f\in \C[\mathfrak{t}]$.

Our previous results, imply $\NH^P\cong e_P \NH e_P $ as graded $H_G^*(pt)$-algebras. 
Furthermore, one has $\NH^P=\NH^{W_P\times W_P}$ as graded $H_G^*(pt)$-modules and  $\NH^{W\times W}=\C[\mathfrak{t}]^W=H_G^*(pt)$. 

\begin{lemma} \label{nilHeckeCorner}
\begin{itemize}
\item[(1)] Let $s= \# W_P$. We have an isomorphism of $\C[\mathfrak{t}]^W$-algebras 
\[\NH\cong M_s(\NH^P),\] 
in particular it is a free module over $\NH^P$ of rank $s^2$. 
\item[(2)] Let $r=\# W^P$. We have an isomorphism of $\C[\mathfrak{t}]^W$-algebras 
\[ \NH^P \cong M_r( \C[\mathfrak{t}]^W ), \]
in particular it is a free module over $\C[\mathfrak{t}]^W$ of rank $r^2$. \\
A basis is given by $c_{v,w}, v,w\in W^P$ with $c_{v,w}$ is a lift of 
$[\overline{BvP/P}\times \overline{BwP/P}] \in H_* (G/P\times G/P)$ to $H_*^G(G/P\times G/P)$, i.e. elements in the fibres 
of the forgetful map (which is a surjective ring homomorphism), 
\[
\begin{aligned}
H_*^G(G/P\times G/P) \cong M_r(\C[\mathfrak{t}]^W)  &\surj 
M_r(\C)\cong \End_{\C} (\C[\mathfrak{t}]^{W_P}/I_W) \cong H_*(G/P\times G/P),\\
(f_{i,j})_{i,j} & \mapsto (f_{i,j}(0))_{i,j}
\end{aligned}
\]
where $I_W$ is the ideal generated by the $W$-invariant polynomials of degree $\geq 1$
\end{itemize}
\end{lemma}
\begin{proof}
\begin{itemize}
\item[(1)] Let $s=\# W_P$, it is the rank of $\C[\mathfrak{t}]$ as module over $\C[\mathfrak{t}]^{W_P}$ and therefore 
\[ \NH \cong \End_{\C[\mathfrak{t}]^W}((\C[\mathfrak{t}]^{W_P})^{\oplus s}) \cong M_s(\NH^P) .\]
\item[(2)] $r=\# W^P$ is the rank of $\C[\mathfrak{t}]^{W_P}$ as module over $\C[\mathfrak{t}]^W$ and the dimension as $\C$-vector space of $\C[\mathfrak{t}]^{W_P}/I_W$. The rest follows as in (1).  
\end{itemize}
\end{proof}

\begin{remark} Let $S\subset W$ denote the simple reflections with respect to the Borel $B$ and $n=\rk T$. For $s\in S$ let 
$\al_s \in \C[x_1, \ldots , x_n]$ be a linear polynomial with $s(\al_s)=-\al_s$. 
We denote by $\delta_s \colon \C[x_1,\ldots ,x_n] \to \C[x_1, \ldots , x_n]$ the operator $\delta_s (f) := \frac{s(f)-f}{\al_s}$ called divided difference operator. It is well-known that the affine nil Hecke algebra $\NH=\NH^B$ (with $B$ a Borel group) is isomorphic to the subalgebra of $\End_{\C[x_1,\ldots , x_n]^W} (\C[x_1, \ldots , x_n] )$ generated by $x_1\cdot , \ldots , x_n\cdot ,\delta_s$, $s\in S$.  
But we do not know any description of generators (and relations) for the parabolic nil Hecke algebra $\NH^P$, the analogues of the divided difference operators are misssing. 
\end{remark}

\subsection{On the parabolic analogue of the Springer map}

\begin{definition}
Let $(G,B,T)$ a complex reductive group with Borel subgroup and maximal torus $T$. 
Let $(W,S)$ be the associated Coxeter system and 
$J\subset S$, we set $W_J:= \langle J\rangle \subset W$ and write $P_J = B W_JB$ for the standard parabolic group. Recall, the affine nil Hecke algebra is the graded $\C$-algebra 
$\NH :=  \End_{\C[\mathfrak{t}]^W} (\C[\mathfrak{t}])$. Now, we define a subalgebra $A_J\subset \NH$ to be the subalgebra generated by multiplication with elements in $\C[\mathfrak{t}]$ and for each $s\in S$ a generator 
\[
\si (s):= \begin{cases} \delta_s, &\text{ if } s\in J \\ s-\id , & \text{ if } s\notin J \end{cases}.
\]
We have for $J\subset K\subset S$ that $A_J\subset A_K$ and for arbitrary two subsets 
$J,K\subset S$ one has $A_J\cap A_K =A_{J\cap K}$ and the subalgebra generated by $A_J$ and $A_K$ is $A_{J\cup K}$. Since $A_{\emptyset } = \C [\mathfrak{t} ]\# \C W$, $A_S= \NH$, we see 
these algebras as a lattice of \textbf{interpolations} between the Steinberg algebra $A_{\emptyset}$ corresponding to the Springer map and the nil Hecke algebra $A_S$. 
\end{definition}

We shorten the notation $P:=P_J$. Now, the \textbf{parabolic version of the classical Springer map}, is given by 
\[ 
\pi\colon G\times^P \mathfrak{u}_P \to \mcN, \overline{\left(g,f \right)} \mapsto gf
\]
with $\mathfrak{u}_P=\Lie U_P$ where $U_P\subset P$ is the unipotent radical. The fibres of this map have been studied in \cite{BM}. Let $\mcZ^P$ be the associated Steinberg algebra (i.e. $\mcZ^P= H^G_{[*]}((G\times^P \mathfrak{u}_P)\times_{\mcN} (G\times^P \mathfrak{u}_P))$). The degree zero part of these algebras has already been described in the main result of Douglas and R\"ohrle, 
see \cite{DR}. Our result gives us the following. 

\begin{lemma} In the notation as before with $P:= P_J$ one has: 
\begin{itemize}
\item[(1)] 
The algebra $A_J$ is the Steinberg algebra $\mcZ^B$ associated to the map 
$\pi\colon G\times^B \mathfrak{u}_P \to \mcN, \overline{(g,f)} \mapsto gf.$
\item[(2)] Let $e:= \frac{1}{\# W_J} \sum_{w\in W_J} w \in \NH$ be the idempotent element given by the averaging map from the first section. One has $e\in A_J$ and 
\[ 
eA_Je = \mcZ^P.
\]
Furthermore, one has $eA_Je =A_J\cap \NH^{P_J}$.
\end{itemize}
\end{lemma}

\paragraph{proof:} 
\begin{itemize}
\item[(1)] We have the monomorphism $\mcZ^B \to \NH$ from before. By the main result from \cite{Sa}, we know that the image is the subalgebra generated by the elements $\si(s), s\in S$ as given above and multiplication by polynomials in $\C[\mathfrak{t}]$.  
\item[(2)] Follows from (1) using our main result. For the last statement, we have to see  
$eA_Je= A_J \cap e\NH e $. Clearly, $eA_Je$ is a subset of $A_J \cap e\NH e$. If you take an 
element $ebe=a\in A_J, b\in \NH$, then $ebe=eae \in eA_Je$.  
\hfill $\Box $
\end{itemize}

\subsection{Reineke's Example} 
In the end of \cite{R}, Reineke looked for a description of a parabolic Steinberg algebra, we give the answer which we get from our main result. 
Let $Q$ be the quiver $(1\to 2)$ and let $(d_1,d_2)\in \N_0^{Q_0}$. A directed partition of the Auslander-Reiten quiver of $\C Q$ is given by 
$ I_1:= \{ E_2:=(0\to \C)\}, I_2:=\{ E_{1,2}:= (\C\xrightarrow{\id} \C)\}, I_3:= \{ E_1:= (\C\to 0)\}$, i.e. it is a partition of the vertices of the Auslander-Reiten quiver $\{I_t\}_t$ such that $\EXT^1 (I_t, I_t)=0$ and $\forall t<u  $ $\Hom (I_u, I_t)=0=\EXT^1(I_t,I_u)$. Let $M= E_2^{d_2}\oplus E_{1,2}\oplus E_{1}^{d_1}$. Then, M. Reineke proved that quiver-graded Springer map 
corresponding to the dimension filtration $(0, (d_1-1, 0), (d_1-1, 1), (d_1,1), (d_1,d_2))$ gives a resolution of singularities (i.e. birational projective map) 
for the orbit closure of $M$. Yet, we will consider the even easier dimension filtration $\ddd :=(0, (d_1,1), (d_1,d_2))$. The associated Steinberg variety is 
\[
Z:= \{ (A, L_1, L_2)\in M_{d_2\times d_1}(\C)\times \P^{d_2}(\C ) \times \P^{d_2}(\C )\mid \Bild (A) \subset L_i, i=1,2\}, 
\]
it carries an operation of $\Gd:= \Gl_{d_2}\times \Gl_{d_1}$ via 
\[(A,L_1,L_2)\mapsto (g_2^{-1}Ag_1, g_2^{-1}L_1, g_2^{-1}L_2), \quad (g_2,g_1)\in \Gd.\] 
We want to describe the Steinberg algebra $H_*^{\Gd}(Z)$ with our method.  
We set  
\begin{itemize}
\item[] $\G:= \Gl_d, d:= d_1+d_2,$
\item[] $T:=$ invertible diagonal matrices,
\item[] $\B:=$ invertible upper triangular matrices, 
\item[] $\P:=$ invertible upper block matrices with diagonal block sizes $(1, d-1)$, 
\item[] $\mcU= \Lie \U_{\P}$, where $\U_{\P}$ is the unipotent radical of $\P$,    
\item[] $G:= \Gl_{d_2}\times \Gl_{d_1}$ diagonally embbeded into $\G$, 
\item[] $V= M_{d_2\times d_1}$ embedded into the right upper corner of $\mcG=\mathfrak{g}l_d$. 
\end{itemize}
as usual set $B:= \B \cap G, P:= \P \cap G$, $F:= \mcU \cap V$.  

The algebra $ H_*^G(Z^B), Z^B:= (G\times^B F)\times_V (G\times^B F)$ can by theorem \cite[theorem 2.1]{Sa}, be described as 
the the algebra $1_e*\mcZ* 1_e$ for $\mcZ$ be the Steinberg algebra associated to 
$(\G, \B, \mcU, V)$ and $e\in W\setminus \W$ be the coset of the neutral element. If we set $s_i:= (i,i+1)\in S_d$ and 
\[\delta_i:= \delta_{s_i}\colon \C[ t_1,\ldots , t_d]\to \C[t_1,\ldots , t_d] , \quad f\mapsto \frac{s_i(f)-f} {t_i-t_{i+1}} \]
Then, $H_G^*(Z^B)$ is the subalgebra of 
$\End_{\C[t_1,\ldots , t_d]^{S_{d_2}\times S_{d_1}}}(\C[t_1,\ldots , t_d])$ generated by \[(t_j\cdot), \;1\leq j\leq d, \quad  \delta_i, \; i\in \{ 2,\ldots d_2-2, d_2,\ldots , d-1\},\] 
\[\theta:= \prod_{j=d_1+d_2+1}^d (t_1-t_j) \delta_{1}. \]
Now, Reineke's variety equals $Z=(G\times^PF)\times_V (G\times^PF)$, by the previous section we conclude it is the corner algebra of 
$e_PH_*^G(Z^B)e_P$ where \[e_P\colon \C[t_1,\ldots, t_d]\to \C[t_1,\ldots , t_d], f\mapsto \frac{1}{(d_2-1)! d_1!} \sum_{w\in <s_2,\ldots s_{d_2-2}, s_{d_2}, \ldots , s_{d-1}> }   w(f).\] 

\paragraph{Acknowlegdement:} This article has been part of my phd thesis at the University of Leeds. I want to thank my supervisor Andrew Hubery for detailed discussions and William Crawley-Boevey for his help and advice. Furthermore, I would like to thank Syu Kato for his remark about graded Morita equivalence.  
I had financial support through a 6-months-grant from the University of Leeds and I want to thank the CRC 701 for guest stays at the Universit\" at Bielefeld.

%
%


\bibliographystyle{alphadin}
\bibliography{IdemRecForParabolics}
%

\end{document}